\documentclass[preprint,1p,times]{elsarticle}
\usepackage{amsmath, amssymb, amsthm}
\usepackage{graphicx}
\usepackage{color,soul}
\usepackage{float}
\usepackage{lineno}
\usepackage{pgf}
\usepackage{epstopdf}

\journal{ }
\begin{document}
	
	\begin{frontmatter}
		
		\title{Analysis of a fractional order eco-epidemiological model with prey infection and type II functional response}
		\author[cmbe]{Shuvojit Mondal\fnref{fn1}}
		\ead{shuvojitmondal91@gmail.com}
		\cortext[cor1]{Corresponding author}
		\author[ju]{Abhijit Lahiri\fnref{fn2}}
		\ead{lahiriabhijit2000@yahoo.com}
		\author[cmbe]{Nandadulal Bairagi\corref{cor1}}
		\ead{nbairagi.math@jadavpuruniversity.in}
		\address[cmbe]{Centre for Mathematical Biology and Ecology \\ Department of Mathematics, Jadavpur University\\ Kolkata-700032, India.}
		\address[ju]{Department of Mathematics,Jadavpur University\\ Kolkata-700032, India.}
		\fntext[cor1]{Research of N. Bairagi is supported by DST PURSE, Phase II.}

		\begin{abstract}
			In this paper, we introduce fractional order into an ecoepidemiological model, where predator consumes disproportionately large number of infected preys following type II response function. We prove different mathematical results like existence, uniqueness, non-negativity and boundedness of the solutions of fractional order system. We also prove the local and global stability of different equilibrium points of the system. The results are illustrated with several examples.
		\end{abstract}
		
		
	\end{frontmatter}

\section{Introduction}
Fractional calculus, which is a generalization of integer order differentiation and $n-$fold integration, has been successfully applied in different branches of science and engineering. Differential equations with fractional-order derivatives (or integrals) are generally called fractional differential (or integral) equations. In recent past, fractional order differential equations have been used in several biological systems to explore the underlying dynamics  \cite{CuiYang14, DasGupta09, DasGupta11}. Here we assume an ecological system where a prey population grows logistically and a predator population feeds on this prey population. Now assume that the prey population is infected by some microparasites. In presence of infection, our prey population is divided into two subpopulations, viz. susceptible prey and infected prey. Since infected preys are weaken and cannot easily escape predation, predators disproportionately consume large number of infected prey \cite{LM96}. In such case, if there are sufficient numbers of infected prey, it may be assumed that the growth rate of predator is maintained mainly by consuming infected prey. If it is also considered that healthy preys can only give birth, infection transmits horizontally and predation process follows type II response function then we have the following eco-epidemiological system:
\begin{eqnarray}\label{Eco-epidemiological integer order model}
\frac{dS}{dt} & = & rS\bigg(1-\frac{S+I}{K}\bigg) - \lambda IS,\nonumber \\
\frac{dI}{dt}& = & \lambda IS-\frac{mIP}{a+I} - \mu I, \\
\frac{dP}{dt}& = & \frac{\theta IP}{a+I}-dP. \nonumber
\end{eqnarray}
The state variables $S(t), I(t)$ and $P(t)$ represent, respectively, the densities of susceptible, infected and predator populations at time $t$. Here $r$ is intrinsic birth rate of prey, $K$ is the environmental carrying capacity, $\lambda$ is the force of infection, $m$ is the maximum prey attack rate, $\mu$ is the death rate of infected prey, $\theta$ ($0<\theta\leq1$) is the conversion efficiency, $a$ is the half saturation constant and $d$ is the death rate of predator, All parameters are assumed to be positive from biological point of view. Readers are referred to \cite{ChattoBairagi01} for more discussion about the model. Note that it is an integer order system of differential equations and its dynamics was studied by Chattopadhyay and Bairagi \cite{ChattoBairagi01}.


Considering the fractional derivatives in the sense of Caputo derivative, and assuming $0<\alpha\leq1$, we have the following fractional order eco-epidemiological model corresponding to the model (\ref{Eco-epidemiological integer order model}):
\begin{eqnarray}\label{Eco-epidemiological fractional order model}
^{c}_{0} D^{\alpha}_{t}S & = & rS\bigg(1-\frac{S+I}{K}\bigg) - \lambda IS, \nonumber \\
^{c}_{0} D^{\alpha}_{t}I & = & \lambda IS-\frac{mIP}{a+I} - \mu I, \\
^{c}_{0} D^{\alpha}_{t}P & = & \frac{\theta IP}{a+I}-dP, \nonumber
\end{eqnarray}
where $^{c}_{0} D^{\alpha}_{t}$ is the Caputo fractional derivative. The main advantage of Caputo's approach is that the initial conditions for the fractional differential equations with Caputo derivatives takes the similar form as for integer-order differential equations \cite{CuiYang14, Podlubny99}, i.e., it has advantage of defining integer order initial conditions for fractional order differential equations. We analyze system (\ref{Eco-epidemiological fractional order model}) with the initial conditions $$S(0)>0, I(0)>0, P(0)>0.$$

The paper is organized as follows. In Section $2$, we give some useful theorems and lemmas in relation to fractional order differential equations. Well-posedness and dynamical behavior of the model are presented in Sections $3$ and $4$, respectively. Extensive numerical computations are presented in Section $5$ and the paper ends with a summary in Section $6$.
\section{Important Results}
\noindent\textbf{Theorem 1.} \cite{Matignon96} The following autonomous system
\begin{equation}\nonumber
^{c}_{t_{0}} D^{\alpha}_{t} x(t) = Ax, ~x(0)= x_{0},
\end{equation}
with $0<\alpha\leq1, x\in \Re^{n}$ and $A\in \Re^{n\times n}$ is asymptotically stable if and only if $\mid arg(\lambda)\mid >\frac{\alpha \pi}{2}$ is satisfied for all eigenvalues of the matrix $A$. Also, this system is stable if and only if $\mid arg(\lambda)\mid \geq \frac{\alpha \pi}{2}$ for all eigenvalues of the matrix $A$ with those critical eigenvalues satisfying $\mid arg(\lambda)\mid =\frac{\alpha \pi}{2}$ having geometric multiplicity of one. The geometric multiplicity of an eigenvalue $\lambda$ of the matrix $A$ is the dimension of the subspace of vectors $v$ for which $Av = \lambda v$.\\

\noindent\textbf{Theorem 2.} \cite{Petras11} Consider the following commensurated fractional order system
\begin{equation}\label{Stability condition} \nonumber
^{c}_{t_{0}} D^{\alpha}_{t} x(t) = f(x), x(0)= x_{0}
\end{equation}
with $0<\alpha\leq1, x\in \Re^{n}$ and $f: \Re^{n}\rightarrow \Re^{n}$ i.e., $f=(f_1, f_2 ----f_n)^T$. The equilibrium points of the  above system are calculated by solving the equation $f(x) = 0$. These equilibrium points are locally asymptotically stable if all eigenvalues $\lambda_{i}$ of the jacobian matrix $J = \frac{\partial f}{\partial x}$ evaluated at the equilibrium points satisfy $\mid arg(\lambda_{i})\mid >\frac{\alpha \pi}{2}$.\\

\noindent\textbf{Lemma 1} \cite{OdibatShawagfeh07}~~~(Generalized Mean Value Theorem) Suppose that$f(t)\in C[a,b]$ and $D^{\alpha}_{a}f(t)\in C(a,b]$ with $0<\alpha\leq1$, then we have
\begin{equation}\nonumber
f(t) = f(a) + \frac{1}{\Gamma(\alpha)}(D^{\alpha}_{a}f)(\xi). (t-a)^\alpha,
\end{equation}
where $a\leq\xi\leq x$, $\forall x \in (a,b]$.\\

\noindent\textbf{Corollary 1}~~~ Suppose $f(t)\in C[a,b]$ and $^{c}_{t_0} D^{\alpha}_{t}f(t)\in C(a,b)$, $0<\alpha \leq 1$. If $^{c}_{t_0} D^{\alpha}_{t}f(t)\geq 0, \forall t\in (a,b)$, then $f(t)$ is a non-decreasing function for all $t\in[a,b]$; and if $^{c}_{t_0} D^{\alpha}_{t}f(t)\leq 0, \forall t\in (a,b)$, then $f(t)$ is a non-increasing function for all $t\in[a,b]$. \\

\noindent\textbf{Lemma 2} \cite{Kilbas06}~~~The solution to the cauchy problem
$$^{c}_{t_0} D^{\alpha}_{t}x(t)  = \lambda x(t) + f(t),$$
$$x(a)                         =  b (b\in \Re)$$
with $0< \alpha \leq 1$ and $\lambda\in\Re$ has the form
\begin{equation}\nonumber
x(t) = bE_{\alpha}[\lambda(t-a)^{\alpha}] + \int_a^t (t-s)^{\alpha-1} E_{\alpha,\alpha}[\lambda(t-s)^{\alpha}]f(s)ds,
\end{equation}
while the solution to the problem
$$^{c}_{t_0} D^{\alpha}_{t}x(t)  = \lambda x(t),$$
$$x(a)                         =  b ~(b\in \Re)$$
is given by
\begin{equation}\nonumber
x(t) = bE_{\alpha}[\lambda(t-a)^{\alpha}].
\end{equation}


\noindent\textbf{Lemma 3} \cite{Hong-Li16}~~~Let $u(t)$ be a continuous function on $[t_{0},\infty)$ and satisfying
$$^{c}_{t_0} D^{\alpha}_{t}u(t) \leq -\lambda u(t) + \mu,$$
$$u(t_{0})                         =  u_{t_{0}} $$
where $0< \alpha<1$, $(\lambda,\mu)\in\Re^{2}$, $\lambda\neq0$ and $t_{0}\geq0$ is the initial time. Then its solution has the form
\begin{equation}\nonumber
u(t) = \bigg(u_{t_{0}}- \frac{\mu}{\lambda}\bigg)E_{\alpha}[-\lambda(t-t_{0})^{\alpha}] + \frac{\mu}{\lambda}.
\end{equation}

\noindent\textbf{Lemma 4} \cite{LiChen10}~~~Consider the system
\begin{equation}\nonumber
^{c}_{t_{0}} D^{\alpha}_{t} x(t) = f(t,x), t>t_{0}
\end{equation}
with initial condition $x_{t_{0}}$, where $0<\alpha\leq1$, $f:[t_{0},\infty)\times\Omega\rightarrow\Re^{n}$, $\Omega\in\Re^{n}$. If $f(t,x)$ satisfies the locally Lipschitz condition with respect to $x$, then there exists a unique solution of the above system on $[t_{0},\infty)\times\Omega$.
\section{Well-posedness}
\subsection{Non-negativity and boundedness}
Considering biological significance of the problem, we are only interested in solutions that are non-negative and bounded. Denote $\Re^3{}_{+} = \{x\in\Re^{3}|x\geq0 \}$ and $x(t) = (S(t),I(t),P(t))^{T}$.\\

\noindent\textbf{Theorem 3.}~~All solutions of the system (\ref{Eco-epidemiological fractional order model}) which start in $\Re^3{}_{+}$ are non negative and uniformly bounded.\\

\noindent\textbf{Proof.} First we prove that $S(t)\geq0, \forall t\geq0$ assuming $S(0)>0$ for $t=0$. Let us suppose that $S(t)\geq0, \forall t\geq0$ is not true. Then there exists some $t_{1}> 0$ such that
$S(t)> 0$ for $0\leq t< t_{1}$, $S(t) = 0$ at $t = t_{1}$ and $S(t)< 0$ for $t> t_{1}$.\\
From the first equation of (\ref{Eco-epidemiological fractional order model}), we have
\begin{equation}
^{c}_{0} D^{\alpha}_{t_{1}} S(t)|_{t=t_{1}} = 0.
\end{equation}
According to Lemma $1$, we the have $S(t_{1}^{+}) = 0$, which contradicts the fact $S(t_{1}^{+}) < 0$, i.e. $S(t)< 0$ for $t> t_{1}$. Therefore, we have $S(t)\geq0, \forall t\geq0$. Using similar arguments, we can prove $I(t)\geq0, \forall t\geq0$ and $P(t)\geq0, \forall t\geq0$. Next we show that all solutions of system (\ref{Eco-epidemiological fractional order model}) which start in $\Re^3{}_{+}$ are uniformly bounded.\\ We define a function
\begin{equation}
V(t) = S + I + \frac{m}{\theta}P,
\end{equation}
Taking its fractional time derivative, we have
\begin{equation}\nonumber
^{c}_{0} D^{\alpha}_{t} V(t) = rS\bigg(1-\frac{S+I}{K}\bigg) - \mu I - \frac{md}{\theta} P.
\end{equation}
Now for each $\eta>0$, we have
\begin{equation}\label{Boundedness}
\begin{split}
^{c}_{0} D^{\alpha}_{t} V(t) + \eta V(t) = & rS\bigg(1-\frac{S+I}{K}\bigg) - \mu I - \frac{md}{\theta} P + \eta S +\eta I + \frac{m\eta}{\theta} P \\
\leq & -\frac{r}{K}S^{2} + (r + \eta)S + (\eta - \mu) I+ (\eta - d)\frac{md}{\theta}P \\
\leq & \frac{K}{4r}(r + \eta)^{2} + (\eta - \mu) I+ (\eta - d)\frac{md}{\theta}P.
\end{split}
\end{equation}
Taking $\eta < min(\mu, d)$, we have
\begin{equation}
\begin{split}
^{c}_{0} D^{\alpha}_{t} V(t) + \eta V(t) \leq l,
\end{split}
\end{equation}
where $l=\frac{K}{4r}(r + \eta)^{2}>0.$ Applying Lemma $3$, one gets
\begin{equation}
\begin{split}
V(t) \leq & (V(0) - \frac{l}{\eta})E_{\alpha}[-\eta t^{\alpha}] + \frac{l}{\eta}, \\
\leq & V(0) E_{\alpha}[-\eta t^{\alpha}] + \frac{l}{\eta}(1 - E_{\alpha}[-\eta t^{\alpha}]).
\end{split}
\end{equation}
Thus,$V(t)\rightarrow \frac{l}{\eta}$ as $t\rightarrow \infty$ and $0<V(t)\leq \frac{l}{\eta}$. Hence all solutions of system (\ref{Eco-epidemiological fractional order model}) that starts from $\Re^3{}_{+}$ are confined in the region $\Omega = \{(S,I,P)\in\Re^3{}_{+}| ~V(t)\leq \frac{l}{\eta} + \epsilon$, for any $\epsilon>0\}$.
\subsection{Existence and uniqueness}
Now, we  study the existence and uniqueness of the solution of system (\ref{Eco-epidemiological fractional order model}) in the region $\Omega \times[0, T]$, where $\Omega = \{(S,I,P)\in\Re^{3}|~ max\{|S|, |I|, |P|\} \leq M\}$, $T<\infty$ and $M$ is sufficiently large. Denote $X = (S,I,P)$, $\bar{X} = (\bar{S}, \bar{I}, \bar{P})$. Consider a mapping $H(X) = (H_{1}(X), H_{2}(X), H_{3}(X))$ and
\begin{eqnarray}\label{Existence}
H_{1}(X) & = & rS\bigg(1-\frac{S+I}{K}\bigg) - \lambda IS, \nonumber \\
H_{2}(X) & = & \lambda IS-\frac{mIP}{a+I} - \mu I, \\
H_{3}(X) & = & \frac{\theta IP}{a+I}-dP. \nonumber
\end{eqnarray}
For any $X, \bar{X} \in \Omega$, it follows from (\ref{Existence}) that
\begin{equation}\nonumber
\begin{split}
\parallel H(X) - H(\bar{X})\parallel =& \mid H_{1}(X)-H_{1}(\bar{X})\mid +\mid H_{2}(X)-H_{2}(\bar{X})\mid +\mid H_{3}(X)-H_{3}(\bar{X})\mid \\
=& \mid rS\bigg(1-\frac{S+I}{K}\bigg) - \lambda IS - r\bar{S}\bigg(1-\frac{\bar{S}+\bar{I}}{K}\bigg) + \lambda \bar{I}\bar{S}\mid \\
& +\mid \lambda IS-\frac{mIP}{a+I} - \mu I - \lambda \bar{I}\bar{S}+\frac{m\bar{I}\bar{P}}{a+\bar{I}} + \mu \bar{I}\mid \\
& +\mid \frac{\theta IP}{a+I}-dP-\frac{\theta \bar{I}\bar{P}}{a+\bar{I}}+d\bar{P}\mid \\
=& \mid r(S-\bar{S}) - \lambda(IS -\bar{I}\bar{S}) - \frac{r}{K}(S^{2} -\bar{S}^{2} +IS - \bar{I}\bar{S})\mid \\
& +\mid \lambda(IS - \bar{I}\bar{S}) - \mu(I - \bar{I}) - m\bigg(\frac{IP}{a+I} - \frac{\bar{I}\bar{P}}{a+\bar{I}}\bigg)\mid \\
& + \mid \theta\bigg(\frac{IP}{a+I} - \frac{\bar{I}\bar{P}}{a+\bar{I}}\bigg) - d(P - \bar{P})\mid \\
\leq& r\mid S-\bar{S}\mid + \lambda\mid IS -\bar{I}\bar{S}\mid + \frac{r}{K}\mid S^{2} -\bar{S}^{2}\mid +\mid IS - \bar{I}\bar{S}\mid \\
& + \lambda\mid IS - \bar{I}\bar{S}\mid + \mu \mid I - \bar{I}\mid + m\mid\frac{IP}{a+I} - \frac{\bar{I}\bar{P}}{a+\bar{I}}\mid \\
& + \theta \mid\frac{IP}{a+I} - \frac{\bar{I}\bar{P}}{a+\bar{I}}\mid + d\mid P - \bar{P}\mid \\
\leq& \bigg(r +\frac{2rM}{K} + (2\lambda + \frac{r}{K})M\bigg)\mid S - \bar{S}\mid \\
& + \bigg((2\lambda + \frac{r}{K})M + \mu + \frac{aM(m+\theta)}{(a+M)^2}\bigg)\mid I-\bar{I}\mid \\
& + \bigg(\frac{aM(m+\theta)}{(a+M)^2}+d +\frac{M^{2}(m+\theta)}{(a+M)^2} \bigg)\mid P-\bar{P}\mid \\
\leq& L\parallel (S,I,P) - (\bar{S}, \bar{I}, \bar{P})\parallel \\
\leq& L\parallel X - \bar{X} \parallel,
\end{split}
\end{equation}
where $L = max\{r +\frac{2rM}{K} + (2\lambda + \frac{r}{K})M,(2\lambda + \frac{r}{K})M + \mu + \frac{aM(m+\theta)}{(a+M)^2}, \frac{aM(m+\theta)}{(a+M)^2}+d +\frac{M^{2}(m+\theta)}{(a+M)^2}\}$.\\ Thus $H(X)$ satisfies Lipschitz condition with respect to $X$ and it follows from Lemma $4$ that there exists a unique solution $X(t)$ of the system (\ref{Eco-epidemiological fractional order model}) with the initial condition $X(0) = (S(0), I(0), P(0))$.
\section{Dynamical behavior}
To obtain equilibrium points of (\ref{Eco-epidemiological fractional order model}), we solve the following simultaneous equations:
\begin{eqnarray}\label{Equilibrium system}
^{c}_{0} D^{\alpha}_{t}S & = & 0\nonumber \\
^{c}_{0} D^{\alpha}_{t}I & = & 0, \\
^{c}_{0} D^{\alpha}_{t}P & = & 0. \nonumber
\end{eqnarray}
We thus obtain $E_{0} = (0,0,0)$ as the trivial equilibrium, $E_{1} = (K,0,0)$ as the axial equilibrium, $E_{2} = (S_1, I_1, 0)$ as the planar equilibrium, where $S_1 = \frac{\mu}{\lambda}$ and $I_1 = \frac{r(\lambda K-\mu)}{\lambda(r + \lambda K)}$ and $E^{*} = (S^{*},I^{*},P^{*})$ as the interior equilibrium, where
\begin{eqnarray}\label{Equilibrium relation}
S^{*} =  K-\bigg(1 + \frac{\lambda K}{r}\bigg)I^{*},~~ I^{*} =  \frac{ad}{\theta -d},~~ P^{*} = \frac{(a+I^{*})(\lambda S^{*} - \mu)}{m}.
\end{eqnarray}
Note that the equilibria $E_{0}$ and $E_{1}$ always exist. The planar equilibrium point $E_{2}$ exists if $R_{0}>1$, where $R_{0} =\frac{\lambda K}{\mu}$. The interior equilibrium $E^{*}$ exists if $(i)$ $R_{0}>1$ and $(ii)$ $\theta > \theta_{1}$, where $\theta_1=d + \frac{\lambda ad(r + \lambda K)}{r(\lambda K - \mu)}$.\\
The jacobian matrix of system (\ref{Eco-epidemiological fractional order model}) evaluated at $E_{0}$ is given by
\[
J(E_{0}) = \begin{pmatrix}
r & 0 & 0 \\ 0 & -\mu & 0 \\ 0 & 0 & -d \end{pmatrix}.\]
The eigenvalues can be determined by solving the characteristic equation $det(J(E_{0} - \xi I_{3})) = 0$ and they are $\xi_{1} = r ~(>0)$, $\xi_{2} = -\mu ~(<0)$ and $\xi_{3} = -d ~(<0)$. Note that $\mid arg(\xi_{1})\mid = 0$, $\mid arg(\xi_{2})\mid = \pi $, and $\mid arg(\xi_{3})\mid = \pi $. Since the first eigenvalue $\xi_{1}$ does not satisfy $\mid arg(\xi_{1})\mid >\frac{\alpha \pi}{2}$ for all $\alpha\in (0,1]$, therefore $E_{0} = (0, 0, 0)$ is always unstable.\\

\noindent The jacobian matrix $J(E_{1})$ is computed as
\[
J(E_{1}) = \begin{pmatrix}
-r & -r-\lambda K & 0 \\ 0 & \lambda K - \mu & 0 \\ 0 & 0 & -d \end{pmatrix}.\]
The corresponding eigenvalues are $\xi_{1} = -r ~(<0)$, $\xi_{2} = \lambda K - \mu$, $\xi_{3} = -d<0$. Here two cases arise depending on whether $R_0>1$ or $R_0<1$. \\
\noindent\textbf{Case 1:} If $R_0<1$ then we can see that $\mid arg(\xi_{i})\mid  = \pi > \frac{\alpha \pi}{2}, \forall \alpha\in (0,1]$, $i = 1, 2, 3$. Therefore, the equilibrium $E_{1}$ is locally asymptotically stable.\\
\noindent\textbf{Case 2:} If $R_0>1$ then it is easy to see that $\mid arg(\xi_{2})\mid = 0$. In this case, $E_{1}$ is unstable. \\

Performing similar calculations, one can show that the characteristic equation of the Jacobian matrix $J(E_{2})$ can be expressed as
\begin{equation}
\bigg(\xi- (\frac{\theta r (\lambda K - \mu)}{a\lambda(\lambda K + r)+ r(\lambda K-\mu)} - d)\bigg)( \xi^{2} + A\xi + B)=0,
\end{equation}
where $A = \frac{r\mu}{\lambda K} > 0$ and $B = \frac{r\mu (\lambda K- \mu)}{\lambda K} > 0$.
Therefore, one eigenvalue is $\xi_{1} = d_{1} - d$, where $d_{1} = \frac{\theta r (\lambda K - \mu)}{a\lambda(\lambda K + r)+ r(\lambda K-\mu)}$ and the other two are given by $\xi_{2, 3} = \frac{1}{2}(- A \pm \sqrt{A^{2} - 4B})$. Following two cases may arise. \\

\noindent\textbf{Case $d > d_1$:} If $d > d_{1}$ then the equilibrium $E_{2}$ is locally asymptotically stable. Again if $1<R_0<1+\frac{r}{4}$ then $A^2 - 4B>0$. In this case, $\xi_{1}$, $\xi_{2}$, $\xi_{3}$ all are real negative and $\mid arg(\xi_{i})\mid  = \pi > \frac{\alpha \pi}{2}, \forall \alpha\in (0,1]$, $i = 1, 2, 3$. Therefore, the equilibrium $E_{2}$ is stable node if $d > d_{1}$ and $1<R_0<1+\frac{r}{4}$. However, $R_0>1+\frac{r}{4}$ gives $A^2 - 4B<0$. Then $\xi_{1}<0$ and $\xi_{2}, \xi_{3}$ become complex conjugate with negative real parts. Thus, $\mid arg(\xi_{1})\mid  = \pi > \frac{\alpha \pi}{2}$ and $\mid arg(\xi_{2,3})\mid  = \arctan(\frac{\sqrt{-(A^{2} - 4B)}}{A})> \frac{\alpha \pi}{2}$, $\forall \alpha\in (0,1]$. Therefore, the equilibrium $E_{2}$ is stable focus if $d > d_{1}$ and $R_0>1+\frac{r}{4}$.\\
\noindent\textbf{Case $d < d_{1}$:} If $d < d_{1}$ then $E_{2}$ is always unstable. It will be unstable node if $1<R_0<1+\frac{r}{4}$ and unstable focus if $R_0>1+\frac{r}{4}$.\\

\noindent For the interior equilibrium $E^{*}$, the Jacobian matrix is evaluated as
\[\mathbf{J(E^{*})} = \begin{pmatrix}
-\frac{r S^{*}}{ K} & -(\frac{r}{K}+\lambda)S^{*} & 0 \\ \lambda I^{*} & \frac{mI^{*}P^{*}}{(a + I^{*})^{2}} & -\frac{md}{\theta} \\ 0 & \frac{a \theta P^{*}}{(a + I^{*})^{2}} & 0
\end{pmatrix}.\]
The eigenvalues are the roots of the cubic equation
\begin{equation}\label{Cubic equation}
F(\xi) = \xi^{3} + A_{1} \xi^{2} + A_{2} \xi + A_{3} = 0,
\end{equation}
where $A_{1} = \frac{r S^{*}}{K} - \frac{mI^{*}P^{*}}{(a + I^{*})^{2}},~
A_{2} = \frac{amdP^{*}}{(a+I^{*})^{2}} + \frac{r \lambda I^{*}S^{*}}{K} + \lambda^{2} I^{*}S^{*} - \frac{rmS^{*}I^{*}P^{*}}{K(a+I^{*})^{2}},~
A_{3} = \frac{rmdaS^{*}P^{*}}{K(a+I^{*})^{2}}.$

\noindent The discriminant $D(F)$ of the cubic polynomial $F(\xi)$ is
\[
\mathbf{D(F)} = - \begin{vmatrix}
1 & A_{1} & A_{2} & A_{3} & 0 \\ 0 & 1 & A_{1} & A_{2} & A_{3} \\ 3 & 2A_{1} & A_{2} & 0 & 0 \\ 0 & 3 & 2A_{1} & A_{2} & 0 \\ 0 & 0 & 3 & 2A_{1} & A_{2}
\end{vmatrix}.
\]
On expansion, one gets $D(F) = 18A_{1}A_{2}A_{3} + (A_{1}A_{2})^2 - 4A_{3}A_{1}^{3} - 4A_{2}^{3} - 27A_{3}^{2}$. We have the following proposition. \\
\noindent\textbf{Proposition 1.}\begin{itemize}
	\item[(i)] If $D(F) > 0$, $ A_{1}>0$, $A_{3}>0$ and $A_{1}A_{2}- A_{3}>0$ then the interior equilibrium $E^{*}$ is locally asymptotically stable.
	\item[(ii)] If $D(F) < 0$, $A_{1}\geq 0$, $A_{2}\geq 0$, $A_{3} > 0$ and $0< \alpha < \frac{2}{3}$ then the interior equilibrium $E^{*}$ is locally asymptotically stable.
	\item[(iii)] If $D(F) < 0$, $A_{1} < 0$, $A_{2} < 0$ and $\alpha > \frac{2}{3}$ then the interior equilibrium $E^{*}$ is unstable.
	\item[(iv)] If $D(F) < 0$, $A_{1} > 0$, $A_{2} > 0$, $A_{1}A_{2} = A_{3}$ and $0< \alpha< 1$ then the interior equilibrium $E^{*}$ is locally asymptotically stable.
\end{itemize}
{\bf Proof.} ~(i) If $D(F)$ is positive then all the roots of (\ref{Cubic equation}) are real and distinct. If not, let us assume that $F(\xi) = 0$ has one real root $\xi_1$ and another two complex conjugate roots $\xi_2$, $\xi_3$. In terms of the roots, the discriminant of $F(\xi)$ can be written as \cite{Janson07}
\begin{equation}\label{Discriminant eqn}
D(F) =[( \xi_1 - \xi_2)(\xi_1 - \xi_3)(\xi_2 - \xi_3)]^2.
\end{equation}
Note that
\begin{eqnarray}
\begin{split}
( \xi_1 - \xi_2)(\xi_1 - \xi_3)(\xi_2 - \xi_3) = & ( \xi_1 - \xi_2)(\xi_1 - \overline{\xi_2})(\xi_2 - \overline{\xi_2})\\
= & ( \xi_1 - \xi_2)(\xi_1 - \overline{\xi_2})2 Im(\xi_2)i\\
= & ( \xi_1 - \xi_2)\overline{(\xi_1 - \xi_2)}2 Im(\xi_2)i\\
= & 2|\xi_1 - \xi_2|^2 Im(\xi_2)i.
\end{split}
\end{eqnarray}
Thus, \begin{equation}
D(F) =[2|\xi_1 - \xi_2|^2 Im(\xi_2)i]^2 < 0,
\end{equation}
which contradicts the fact that $D(F) > 0$. Therefore, whenever $D(F)>0$ then $F(\xi) = 0$ has three real distinct roots. Since $ A_{1}>0$, $A_{3}>0$ and $A_{1}A_{2}- A_{3}>0$, all roots of $F(\xi) = 0$ has negative real roots or complex conjugate roots with negative real parts. As $D(F) > 0$, so all roots of $F(\xi) = 0$ are real negative. Consequently, $\mid arg(\xi_{i})\mid  = \pi > \frac{\alpha \pi}{2}, \forall \alpha\in (0,1]$, $i = 1, 2, 3$, and the equilibrium $E^{*}$ is  locally asymptotically stable. This completes the proof of (i).\\

\noindent (ii) We have seen in (i) that $F(\xi) = 0$ has one real and two complex conjugate roots if $D(F) < 0$. Since $A_{3} > 0$, following (\ref{Cubic equation}), the real root is negative. We thus consider the roots as $\xi_{1} = -b$, $(b\in R_+)$ and $\xi_{2,3} = \beta \pm i \gamma,$ $(\beta, \gamma \in R)$ and
$$F(\xi) = (\xi +b)(\xi - \beta - i \gamma)(\xi - \beta + i \gamma).$$ Comparing this with (\ref{Cubic equation}), we have
$A_{1} =  b-2\beta,
A_{2} = \beta^2 + \gamma^2 - 2b\beta,
A_{3} = b(\beta^2 + \gamma^2).$
Now $A_{1} \geq 0 \Rightarrow b \geq 2\beta$. Noting $\beta^2 \sec^2\theta= \beta^2 + \gamma^2$ and $A_{2} \geq 0$, we have $sec^2\theta \geq  4.$ Therefore, $\theta= |arg(\xi)| \geq \frac{\pi}{3}$. Since $0< \alpha < \frac{2}{3}$, then $ |arg(\xi)| = \theta \geq \frac{\pi}{3} > \frac{\alpha \pi}{2}$ holds. Thus, all roots of (\ref{Cubic equation}) satisfy $\mid arg(\xi_{i})\mid >\frac{\alpha \pi}{2}, \forall \alpha\in (0,1]$ and the equilibrium $E^{*}$ is locally asymptotically stable. This completes the proof of (ii). Proof of $(iii)$ is similar to the proof of (ii) and hence omitted.\\

\noindent Since $D(F) < 0$, $A_{1} > 0$, $A_{2} > 0$, from the previous case, we have the
\begin{eqnarray}\nonumber
A_{1} =  b-2\beta, ~A_{2} =  \beta^2 + \gamma^2 - 2b\beta, ~A_{3} =  b(\beta^2 + \gamma^2).
\end{eqnarray}
Note that $A_{1} > 0 \Rightarrow b> 2\beta$, $A_{2} > 0 \Rightarrow \beta^2 + \gamma^2 - 2b\beta > 0$
and $A_{1}A_{2} = A_{3} \Rightarrow (b-2\beta)(\beta^2 + \gamma^2 - 2b\beta) = b(\beta^2 + \gamma^2) \Rightarrow \beta(b^2 + \beta^2 + \gamma^2 - 2b\beta) = 0$. Then two cases arise: \\
\noindent\textbf{Case 1:} If $\beta = 0$ then three roots $\xi_{1},\xi_{2},\xi_{3}$ of (\ref{Cubic equation}) are $-b, \pm i\gamma$. One can see that $\mid arg(\xi_{1})\mid  = \pi > \frac{\alpha \pi}{2}$ and $\mid arg(\xi_{2,3})\mid  = \frac{\pi}{2} > \frac{\alpha \pi}{2}, \forall\alpha\in (0,1)$. Therefore, the equilibrium $E^{*}$ is locally asymptotically stable.\\
\noindent\textbf{Case 2:} If $b^2 + \beta^2 + \gamma^2 -2b\beta=0$ then we have $b=\beta$ and $\gamma=0$. Using it in $b> 2\beta$ and $\beta^2 + \gamma^2 >2b\beta$, we obtain $b<0$, which contradicts the assumption $b\in R_+$.\\
Thus, if $A_{1} > 0$, $A_{2} > 0$, $A_{1}A_{2} = A_{3}$ then one root is real negative and the other two are purely imaginary and therefore $\mid arg(\xi_{1})\mid  = \pi > \frac{\alpha \pi}{2}$ and $\mid arg(\xi_{2,3})\mid  = \pi/2 > \frac{\alpha \pi}{2},$ $\forall \alpha\in (0,1)$, implying local asymptotic stability of $E^{*}$. This completes the proof.
\section{Global asymptotic stability}
We now prove the global stability of different equilibrium points of the system (\ref{Eco-epidemiological fractional order model}). \\

\noindent\textbf{Lemma 5} \cite{Vargas15}~~~Let $x(t)\in \Re_{+}$ be a continuous and derivable function. Then for any time instant $t>t_{0}$
\begin{equation}\nonumber
^{c}_{t_{0}} D^{\alpha}_{t}\bigg[x(t) - x^{*} - x^{*}ln\frac{x(t)}{x^{*}}\bigg] \leq \bigg(1-\frac{x^{*}}{x(t)}\bigg)~~{^{c}_{t_{0}}} D^{\alpha}_{t}x(t), x^{*}\in \Re_{+}, \forall \alpha\in(0,1).
\end{equation}

\noindent\textbf{Theorem 4:} The axial equilibrium $E_{1}=(K,0,0)$ is global asymptotically stable if $R_{0} < 1$. \\

\noindent {\bf Proof.} We consider the following Lyapunov function
\begin{equation}
V(S, I, P) = \bigg(S - K - Kln\frac{S}{K}\bigg) + I+ \frac{m}{\theta}P.
\end{equation}
Here $V$ is a $C^1$ function such that $V > 0$ for all values of $(S(t), I(t), P(t))\neq (K,0,0)$ and $V = 0$ only at $(K,0,0)$. Calculating the $\alpha$ order derivative of $V(S, I, P)$~ along the solutions of (\ref{Eco-epidemiological fractional order model}) and using the Lemma $4$, we have
\begin{eqnarray}\label{Globality_1}
\begin{aligned}
^c_{0}D^{\alpha}_{t}V(S, I, P) \leq& \frac{(S - K)}{S} ~{^c_{0}}D^{\alpha}_{t}S(t) + {^c_{0}}D^{\alpha}_{t}I(t)+\frac{m}{\theta} ~{^c_{0}}D^{\alpha}_{t}P(t)\\
=& (S - K)[r(1-\frac{S+I}{K}) - \lambda I] + [\lambda IS-\frac{mIP}{a+I} - \mu I]+ \frac{m}{\theta}[\frac{\theta IP}{a+I}-dP] \\
=& (S - K)[r(1-\frac{S+I}{K})] - \lambda I(S - K)+ [\lambda IS-\frac{mIP}{a+I} - \mu I]+ [\frac{m IP}{a+I}-\frac{dm}{\theta}P] \\
\leq & (S - K)[r(1-\frac{S}{K})] - \lambda I(S - K)+ [\lambda IS-\frac{mIP}{a+I} - \mu I]+ [\frac{m IP}{a+I}-\frac{dm}{\theta}P] \\
=& -\frac{r}{K}(S - K)^{2} + (\lambda K - \mu)I - \frac{dm}{\theta}P. \\
\end{aligned}
\end{eqnarray}
Note that if $R_{0} < 1$ then $^c_{0}D^{\alpha}_{t}V(S,I,P) \leq0, \forall (S, I, P)\in R^{3}_{+}$, and $^c_{0}D^{\alpha}_{t}V(S,I,P) = 0$ at $E_1$. Therefore, the only invariant set on which $^c_{0}D^{\alpha}_{t}V(S,I,P) = 0$ is the singleton $\{E_{1}\}$. Then using Lemma $4.6$ in \cite{Huo15}, which generalizes the integer-order LaSalles Invariance Principle to fractional-order system, it follows that every nonnegative solution tends to $E_{1}$ when $R_{0} < 1$. Thus, $E_{1}$ is global asymptotically stable if $R_{0} < 1$. \\

\noindent\textbf{Theorem 5:} The planner equilibrium $E_{2} = (S_1, I_1, 0)$ is global asymptotically stable if $d > d_{2}$, where $d_{2} = \frac{\theta r (\lambda K - \mu)}{a\lambda(r+ \lambda K)}$.\\

\noindent {\bf Proof.} Let us define the Lyapunov function as
\begin{equation}
V(S, I, P) = \bigg(S - S_{1} - S_{1}ln\frac{S}{S_{1}}\bigg) + \bigg(I - I_{1} - I_{1}ln\frac{I}{I_{1}}\bigg)+ \frac{m}{\theta}P.
\end{equation}
Here $V$ is a $C^1$ function such that $V > 0$ for all values of $(S(t), I(t), P(t)) \neq (S_1, I_1, 0)$~ and $V = 0$ only at $(S(t), I(t), P(t))= (S_1, I_1, 0)$. As before, we have
\begin{eqnarray}
\begin{aligned}
^c_{0}D^{\alpha}_{t}V(S, I, P) \leq& \frac{(S - S_{1})}{S}~ {^c_{0}}D^{\alpha}_{t}S(t) + \frac{(I - I_{1})}{I} ~{^c_{0}}D^{\alpha}_{t}I(t) + \frac{m}{\theta} ~{^c_{0}}D^{\alpha}_{t}P(t)\\
=& (S - S_{1})[r(1-\frac{S+I}{K}) - \lambda I] + (I - I_{1})[\lambda S-\frac{mP}{a+I} - \mu]+ \frac{m}{\theta}[\frac{\theta IP}{a+I}-dP] \\
=& (S - S_{1})[r(\frac{S_{1}+I_{1}}{K}-\frac{S+I}{K}) - \lambda (I - I_{1})] + (I - I_{1})[\lambda (S -S_{1})-\frac{mP}{a+I}] + [\frac{mIP}{a+I}-\frac{dm}{\theta}P] \\
=& -\frac{r}{K}(S - S_{1})^{2} - \frac{r}{K}(S - S_{1})(I - I_{1}) + [\frac{m I_{1}}{a +I} - \frac{dm}{\theta}]P \\
\leq & -\frac{r}{K}(S - S_{1})^{2} - \frac{r}{K}[\frac{(S - S_{1})^{2} + (I - I_{1})^{2}}{2}] + [\frac{m I_{1}}{a +I} - \frac{dm}{\theta}]P \\
= & -\frac{3r}{2K}(S - S_{1})^{2} - \frac{r}{2K}(I - I_{1})^{2} + [\frac{m I_{1}}{a +I} - \frac{dm}{\theta}]P \\
\leq & -\frac{3r}{2K}(S - S_{1})^{2} - \frac{r}{2K}(I - I_{1})^{2} + [\frac{m I_{1}}{a} - \frac{dm}{\theta}]P.\\
\end{aligned}
\end{eqnarray}
One can easily show that if $d >d_2$, where $d_2= \frac{\theta r (\lambda K - \mu)}{a\lambda(r+ \lambda K)}$ then $^c_{0}D^{\alpha}_{t}V(S,I,P) \leq0, \forall (S, I, P)\in R^{3}_{+}$, and $^c_{0}D^{\alpha}_{t}V(S,I,P) = 0$ at $E_2$. Therefore, the only invariant set on which $^c_{0}D^{\alpha}_{t}V(S,I,P) = 0$ is the singleton $\{E_{2}\}$. Following Lemma $4.6$ in \cite{Huo15}, it follows that if $E_{2}$ exists and $d > d_{2}$ then it is global asymptotically stable.  \\

\noindent \textbf{Remark:}
It is to be noted that $d_{2} > d_{1}$, where $d_{1} = \frac{\theta r (\lambda K - \mu)}{a\lambda(\lambda K + r)+ r(\lambda K-\mu)}$. This shows that global stability of $E_{2}$ implies its local stability but the converse is not necessarily true. There may exist some parametric space where $E_{2}$ is only locally stable.\\

\noindent\textbf{Theorem 6:} The positive interior equilibrium $E^{*}=(S^*, I^*, P^*)$ is global asymptotically stable if ~$\theta_{1} < \theta < \theta_{2}$, where $\theta_{1} = d + \frac{\lambda ad(r + \lambda K)}{r(\lambda K - \mu)}$ and $\theta_{2} = \frac{mdK}{2K(\lambda S^{*} - \mu) - r}$. \\
{\bf Proof.} To prove global stability of $E^*$, we define the Lyapunov function as
\begin{equation}
V(S, I, P) = \bigg(S - S^{*} - S^{*}ln\frac{S}{S^{*}}\bigg) + \bigg(I - I^{*} - I^{*}ln\frac{I}{I^{*}}\bigg)+ \frac{m}{\theta}\bigg(P - P^{*} - P^{*}ln\frac{P}{P^{*}}\bigg),
\end{equation}
where $V$ is a $C^1$ function such that $V > 0$ for all values of $(S(t), I(t), P(t)) \neq (S^{*},I^{*},P^{*})$ and $V = 0$ only at $(S(t), I(t), P(t)) = (S^{*},I^{*},P^{*})$. We then have
\begin{eqnarray}\nonumber
\begin{aligned}
^c_{0}D^{\alpha}_{t}V(S, I, P) \leq& \frac{(S - S^{*})}{S} ~{^c_{0}}D^{\alpha}_{t}S(t) + \frac{(I - I^{*})}{I} ~{^c_{0}}D^{\alpha}_{t}I(t)+\frac{m}{\theta}\frac{(P - P^{*})}{P} ~{^c_{0}}D^{\alpha}_{t}P(t)\\
=& (S - S^{*})[r(1-\frac{S+I}{K}) - \lambda I] + (I - I^{*})[\lambda S-\frac{mP}{a+I} - \mu]+ \frac{m}{\theta}(P - P^{*})[\frac{\theta I}{a+I}-d] \\
=& (S - S^{*})[r(\frac{S^{*}+I^{*}}{K}-\frac{S+I}{K}) - \lambda (I - I^{*})] + (I - I^{*})[\lambda (S -S^{*})-\frac{mP}{a+I}+ \frac{mP^{*}}{a+I^{*}}]\\
&+(P - P^{*})[\frac{mI}{a+I}-\frac{mI^{*}}{a+I^{*}}]\\
=& -\frac{r}{K}(S - S^{*})^{2} - \frac{r}{K}(S - S^{*})(I - I^{*}) - \frac{m(I - I^{*})(PI^{*} - P^{*}I)}{(a + I)(a + I^{*})}\\
=& -\frac{r}{K}(S - S^{*})^{2} - \frac{r}{K}(S - S^{*})(I - I^{*}) + \frac{mP^{*}(I - I^{*})^{2}}{(a + I)(a + I^{*})}
- \frac{mI^{*}(I - I^{*})(P - P^{*})}{(a + I)(a + I^{*})}\\
\leq & -\frac{r}{K}(S - S^{*})^{2} - \frac{r}{K}[\frac{(S - S^{*})^{2} + (I - I^{*})^{2}}{2}] + \frac{mP^{*}(I - I^{*})^{2}}{(a + I)(a + I^{*})}\\
&- \frac{mI^{*}}{(a + I)(a + I^{*})}[\frac{(I - I^{*})^{2} + (P - P^{*})^{2}}{2}]\\
= & -\frac{3r}{2K}(S - S^{*})^{2} + [\frac{mP^{*}}{(a + I)(a + I^{*})} - \frac{mI^{*}}{2(a + I)(a + I^{*})} - \frac{r}{2K}](I - I^{*})^{2} - \frac{mI^{*}(P - P^{*})^{2}}{2(a + I)(a + I^{*})}\\
\leq & -\frac{3r}{2K}(S - S^{*})^{2} + [\frac{mP^{*}}{(a + I^{*})} - \frac{mI^{*}}{2(a + I^{*})} - \frac{r}{2K}](I - I^{*})^{2} - \frac{mI^{*}}{2(a + I^{*})}(P - P^{*})^{2}.\\
\end{aligned}
\end{eqnarray}
Observe that $[\frac{mP^{*}}{(a + I^{*})} - \frac{mI^{*}}{2(a + I^{*})} - \frac{r}{2K}] < 0$ if $\theta < \theta_{2}$, where $\theta_{2} = \frac{mdK}{2K(\lambda S^{*} - \mu) - r}$. In this case, $^c_{0}D^{\alpha}_{t}V(S,I,P) \leq0, \forall (S, I, P) \in R^{3}_{+}$ and $^c_{0}D^{\alpha}_{t}V(S,I,P) = 0$ at $(S^*, I^*, P^*)$. Therefore, the only invariant set on which $^c_{0}D^{\alpha}_{t}V(S,I,P) = 0$  is $\{E^{*}\}$. Following Lemma $4.6$ in \cite{Huo15}, whenever the interior equilibrium $E^{*}$ exists and $\theta_{1} < \theta < \theta_{2}$, where $\theta_{1} = d + \frac{\lambda ad(r + \lambda K)}{r(\lambda K - \mu)}$, $\theta_{2} = \frac{mdK}{2K(\lambda S^{*} - \mu) - r}$ then it is global asymptotically stable.
\section{Numerical Simulations}
In this section, we perform extensive numerical computations of our system (\ref{Eco-epidemiological fractional order model}) for different fractional orders $0 < \alpha \leq 1$. We employ Adamas-type predictor corrector method for our fractional order differential equation (FODE) \cite{Diethelm02, Diethelm04}. We first
replace the FODE system (\ref{Eco-epidemiological fractional order model}) by the equivalent fractional order integral system
\begin{eqnarray}\label{Eco-epidemiological fractional integral eqn}
S(t) & = & S(0) + D^{-\alpha}_{t} \bigg(rS[1-\frac{S+I}{K}] - \lambda IS \bigg), \nonumber \\
I(t) & = & I(0) + D^{-\alpha}_{t} \bigg(\lambda IS-\frac{mIP}{a+I} - \mu I \bigg), \\
P(t) & = & P(0) + D^{-\alpha}_{t} \bigg(\frac{\theta IP}{a+I}-dP \bigg), \nonumber
\end{eqnarray}
and then apply the PECE (Predict, Evaluate, Correct, Evaluate) method. With the following three examples we substantiate our analytical findings.

\noindent \textbf{Example 1:}
We consider the parameter values as $r = 2.0$, $K = 40.0$, $\lambda = 0.015$, $m = 0.52$, $\mu = 0.28$, $a = 15.0$, $\theta = 0.189$, $d = 0.09$ and initial point $S(0) = 30, I(0) = 5, P(0) = 10$. Most of the parameter values are taken from \cite{ChattoBairagi01}. Step size for all simulations is considered as $0.05$. We compute that $D(F) = 0.0077 > 0$, $ A_{1} = 1.0879 >0$, $A_{3} =  0.0028>0$, $A_{1}A_{2}- A_{3} = 0.2909>0$. Thus, following Proposition 2(i), the interior equilibrium $E^{*}$ is stable for $0 < \alpha \leq 1$. In Fig. 1 we plot the solutions of FODE system (\ref{Eco-epidemiological fractional order model}) with different values of $\alpha$. It shows that all populations remain stable for all values of $\alpha$ though solutions reach to equilibrium value more slowly for smaller value of $\alpha$.
\begin{figure}[H]
	\includegraphics[width=6.5in, height=2in]{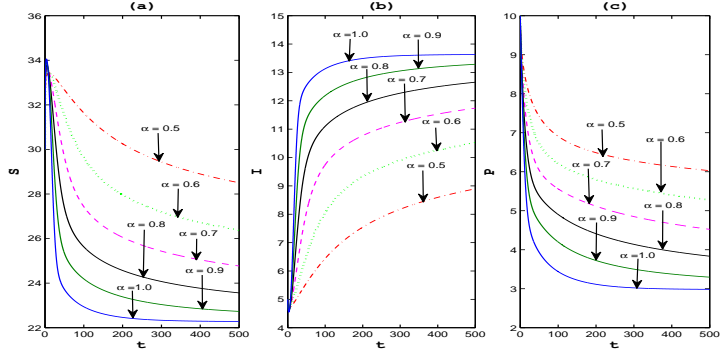}\\
	\caption{Asymptotically stable solutions of $S$, $I$ and $P$ populations with different fractional orders $ 0< \alpha < 1$ and standard order $\alpha = 1$. Here $r = 2.0$, $K = 40.0$, $\lambda = 0.015$, $m = 0.52$, $\mu = 0.28$, $a = 15.0$, $\theta = 0.189$, $d = 0.09$.}
	\label{Stable_1.eps}
\end{figure}
For the above parameter values, following Theorem 6, we determine the critical values of the parameter $\theta$ as $\theta_{1} = 0.1723$ and $\theta_{2} = 0.8044$. Fig. 2 demonstrates that solutions starting from different initial values converge to the equilibrium point $E^{*} = (35.7195, 3.2927, 8.9983)$ for $\theta = 0.5$, depicting its global stability.
\begin{figure}[H]
	\centering
	\includegraphics[width=5in, height=2in]{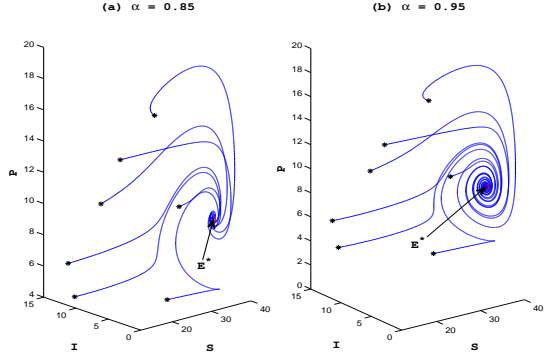}
	\caption{Trajectories with different initial values converge to the coexistence equilibrium $E^{*}$. These figures indicate the global stability of the equilibrium $E^{*}$ for $(a)$ $\alpha = 0.85$ and $(b)$ $\alpha = 0.95$. Here $\theta = 0.5$ and other parameters are as in Fig. 1.}
	\label{Stable_phase_1.eps}
\end{figure}
If we consider $K = 200, \lambda = 0.15, a = 5.0$ and $\theta = 0.9$, keeping other parameter values unchanged, we notice that the conditions of Proposition 2(iii) are satisfied with $D(F) = -463.8995 <0$, $ A_{1} = -0.9276<0$, $A_{2} = -0.5775<0$. Therefore, the interior equilibrium point $E^{*}$ is unstable for  $\alpha > \frac{2}{3}$ (Fig. 3).

\begin{figure}[H]
	\centering
	\includegraphics[width=6.5in, height=2in]{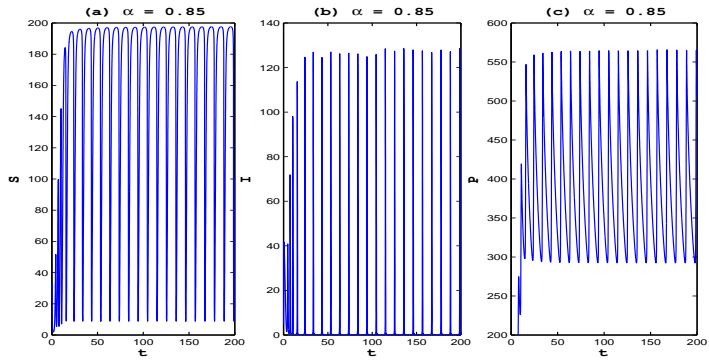}\\
	\caption{Unstable solutions of $S$, $I$ and $P$ populations for $\alpha = 0.85$. Here $K = 200, \lambda = 0.15, a = 5.0$, $\theta = 0.9$ and other parameters are as in Fig. 1.}
	\label{Unstable.eps}
\end{figure}

\noindent \textbf{Example 2.}
For a lower value $\theta= 0.08$, we compute $R_{0} = \frac{\lambda K}{\mu} = 2.142> 1$ and $d - d_1 = 0.0025 > 0$. Therefore, following Theorem 5, the equilibrium point $E_{2} = (18.67, 16.4, 0)$ is stable. Time series solutions and phase portraits of the system (\ref{Eco-epidemiological fractional order model}) for different orders are presented in Figure 4 to illustrate the system behavior. Time evolutions (upper panel) show that solutions converge to the equilibrium faster for higher order and phase diagrams (lower panel) indicate that all trajectories with different initial conditions converge to the predator-free equilibrium $E_{2}$, depicting its is global asymptotic stability.\\

\noindent \textbf{Example 3.}
We now consider the same parameter values and initial point as in Ex. 1 except $\lambda = 0.005$. In this case $R_{0} = \frac{\lambda K}{\mu} = 0.7143< 1$ and we observe that all trajectories with different initial conditions converge to the equilibrium $E_{1} = (40, 0, 0)$ (Fig. 5), following Theorem 4. This indicates that the predator-free and infection-free equilibrium $E_{1}$ is globally asymptotically stable for different orders.

\begin{figure}[H]
	\centering
	\includegraphics[width=6in, height=2.5in]{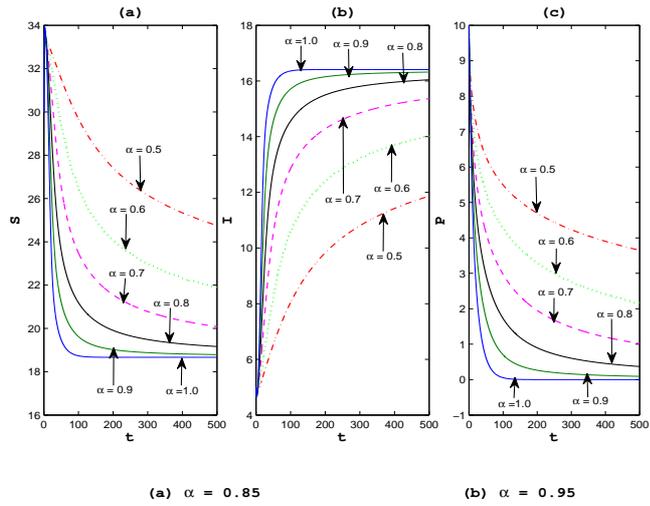}\\
	\includegraphics[width=6in, height=2.5in]{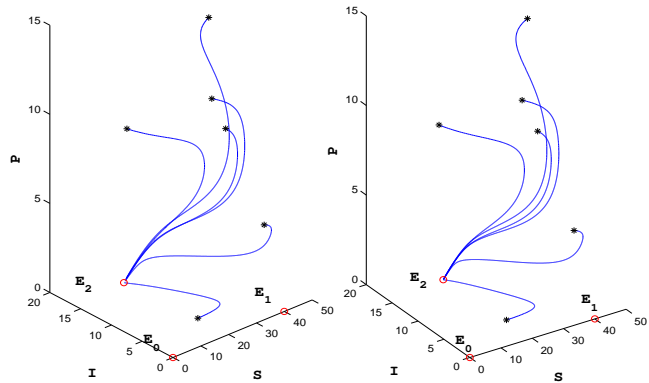}
	\caption{Upper panel shows time behavior and lower panel shows phase behavior of solutions of system (\ref{Eco-epidemiological fractional order model}) for different orders. These figures indicate global stability of predator-free equilibrium $E_{2}$ for different orders. All parameters are as in Fig. 1.}
	\vspace{2in}
	\label{Stable_phase_3.eps}
\end{figure}
\begin{figure}[H]
	\centering
	\includegraphics[width=6in, height=2.5in]{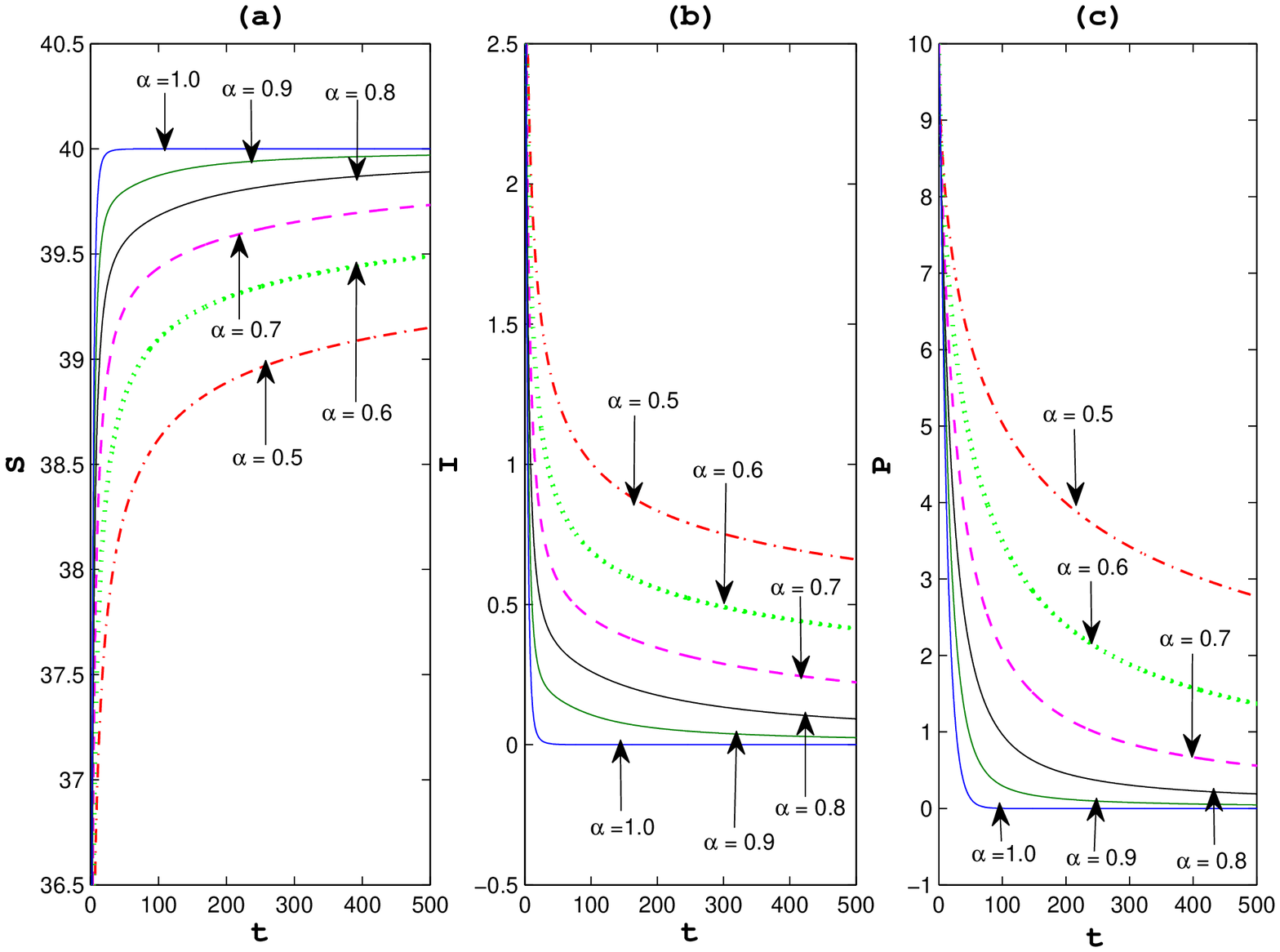}\\
	\includegraphics[width=6in, height=2.5in]{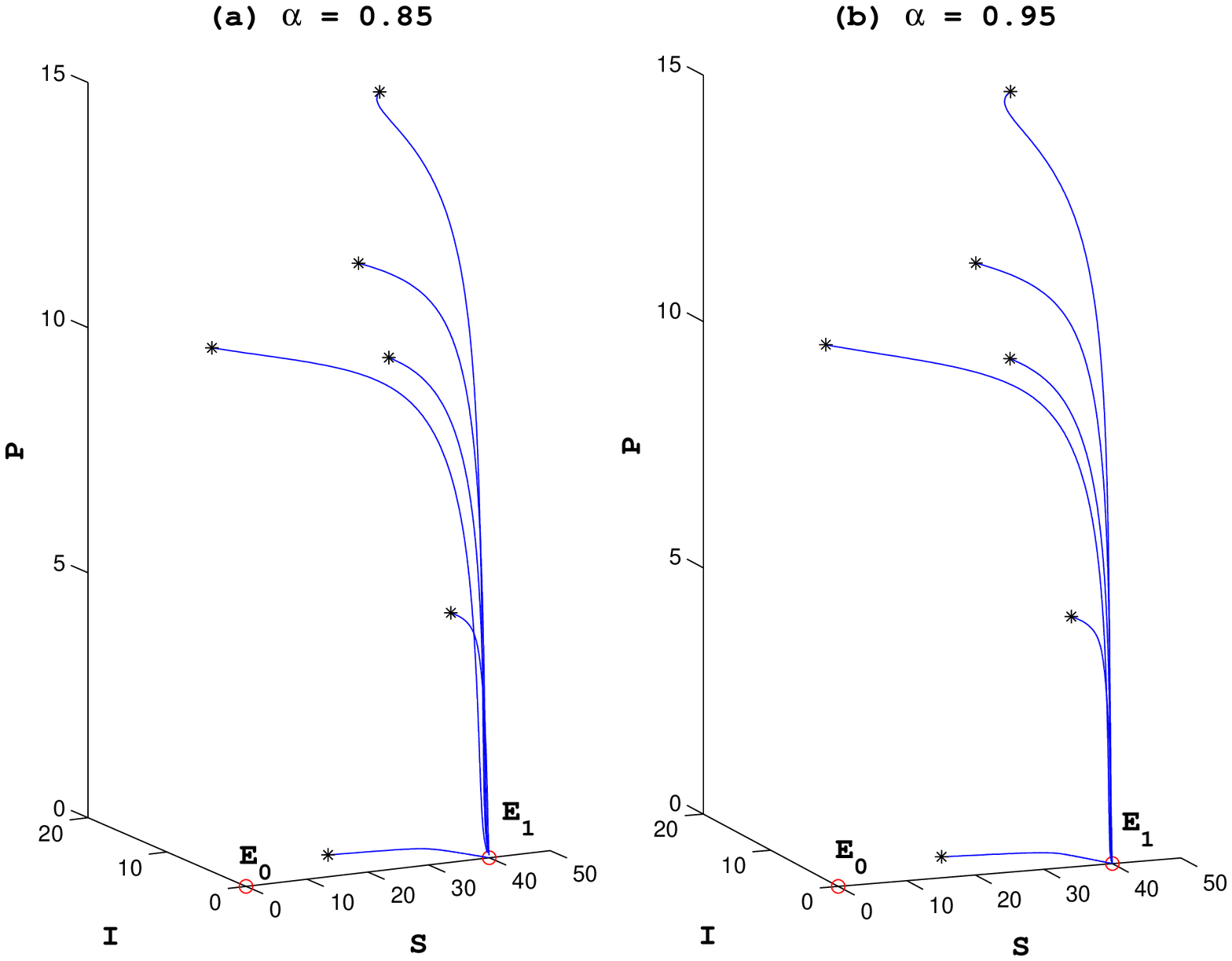}
	\caption{Upper panel shows time behavior and lower panel shows phase behavior of solutions of system (\ref{Eco-epidemiological fractional order model}) for different orders. These figures indicate global stability of the equilibrium $E_{1}$ for different orders. Here $\lambda = 0.005$ with other parameters as in Fig. 1.}
	\label{Stable_phase_2.eps}
\end{figure}

\section{Discussion}
In recent past, eco-epidemiological models have received tremendous attention of modelers because these models consider the issues of ecology and epidemiology simultaneously. Various continuous-time models \cite{Arino99,Ventu94} and discrete-time models \cite{Hu17, Hu14, Adak15} have been proposed and analyzed considering different attributes of the eco-epidemiological system. In this paper, we consider an ecological system where prey population grows logistically and predator population feeds on it following type II response function. When prey is infected by some micro-parasites, it is assumed that predator consumes infected prey only as they are weaken by the disease and can not escape predation. This eco-epidemiological situation has been modeled by a system of fractional order nonlinear differential equations. We prove that the solution of this model system exists uniquely and all solutions remain positive and bounded whenever they start with positive initial value, thus justifying the well-posedness of a biological model. We showed that our system contains four equilibrium points. The trivial equilibrium point is always unstable, implying that all populations can not go to extinction simultaneously. The infection- free and predator-free equilibrium is locally and globally asymptotically stable if $R_0 < 1$ and the dynamics is independent of the order of the differential equation. The predator-free equilibrium is locally asymptotically stable for all order if $1<R_0<1+\frac{r}{4}$ and $d>d_1$. However, if the death rate of predator is very high ($d>d_2>d_1$) then $E_2$ is globally asymptotically stable whenever it exists. The coexistence or interior equilibrium exists if $R_0>1$ and $\theta>\theta_1$. By using stability analysis of fractional order system, we have given different sufficient conditions on the system parameters to prove local stability and instability of $E^*$ for different values of the order, $\alpha$. If, however, $\theta_1<\theta<\theta_2$ then the interior equilibrium is globally asymptotically stable for any $\alpha(0<\alpha\leq1)$ whenever it exists. Numerical examples are presented in support of our analytical results. It is observed that solutions converge to the respective equilibrium more slowly as the order of the differential equation becomes smaller, though the qualitative nature of the solutions remain unchanged.\\

\bibliography{}   

\end{document}